\documentclass[12pt]{article}

\usepackage{amsmath, amsfonts, amssymb}
\usepackage{latexsym}
\usepackage{graphicx}
\usepackage{color}

\newtheorem{thm}{Theorem}[section]
 \newtheorem{cor}[thm]{Corollary}
 \newtheorem{lem}[thm]{Lemma}
 \newtheorem{prop}[thm]{Proposition}
 \newtheorem{exa}[thm]{Example}
 \newtheorem{rem}[thm]{Remark}

\newcommand{\iy}{\infty}

\newcommand{\e}{{\bf e}}
\newcommand{\al}{\alpha}
\newcommand{\be}{\beta}
\newcommand{\ga}{\gamma}

\newcommand{\eps}{\varepsilon}

\newcommand{\ph}{\varphi}

\newcommand{\si}{\sigma}

\newcommand{\vsk}{\vspace{1mm}}
\newcommand{\vsg}{\vspace{2mm}}

\newcommand{\ii}{{\rm i}}
\newcommand{\ee}{{\rm e}}

\renewcommand{\thefootnote}{\fnsymbol{footnote}}

\begin{document}

\begin{center}
{\Large \bf On Hurwitz stable polynomials with integer\\[0.5ex] coefficients}

\vspace{7mm}
{\Large Albrecht B\"ottcher}
\end{center}

%\medskip
%\begin{center}
%{\em To the 90th anniversary of Fekete's subadditive lemma}
%\end{center}

\medskip
\begin{quote}
\footnotesize{Let $H(N)$ denote the set of all polynomials with positive integer coefficients which have their zeros
in the open left half-plane. We are looking for polynomials in $H(N)$ whose largest coefficients are as small as
possible and also for polynomials in $H(N)$ with minimal sum of the coefficients. Let $h(N)$ and $s(N)$ denote
these minimal values. Using Fekete's subadditive lemma we show that the $N$th square roots of $h(N)$ and $s(N)$
have a limit as $N$ goes to infinity and that these two limits coincide. We also derive tight bounds for the
common value of the limits.
}
\let\thefootnote\relax\footnote{\hspace*{-7.5mm} MSC 2010: Primary 12D10; Secondary 12D05, 26C10, 47A68, 93B25}
\let\thefootnote\relax\footnote{\hspace*{-7.5mm} Keywords: Hurwitz polynomial, integer coefficients, Fekete's subadditive lemma}
\end{quote}

\section{Introduction}\label{S1}

A polynomial $p(z)=p_N z^N + \cdots + p_0$ with real coefficients is called Schur stable if all its zeros are in the open unit disk
and is said to be Hurwitz stable if its zeros are all located in the left open half plane. Such polynomials
appear as the result of Wiener-Hopf and spectral factorizations. To test numerical algorithms for these factorizations,
it is desirable to have some supply of Schur and Hurwitz stable polynomials. For example, one starts with
a Hurwitz stable polynomial $p(z)$, forms the product $f(z)=p(-z)p(z)$, applies the algorithm to get a factorization
$f(z)\approx q(-z)q(z)$, and finally one measures the error $q(z)-p(z)$.

\vsk
It is easy to produce nice Schur stable polynomials of arbitrary degree. For instance, Bini, Fiorentino, Gemignani,
and Meini \cite{BiniEff} introduced
the beautiful polynomials
\[p(z)=1+z+\cdots+z^{N-1}+2z^N.\]
To reveal that $p(z)$ is Schur stable, we show that the reverse polynomial
\[z^N p(1/z)=2+z+ \cdots + z^N\] has no zeros in the open unit disk. And indeed, for $z \neq 1$ we have
\[2+z+ \cdots + z^N=1+\frac{1-z^{N+1}}{1-z}=\frac{2-z-z^N}{1-z},\]
and this cannot be zero for $|z|<1$ because then $|z+z^N| < 2$. In fact, the zeros of $p(z)$ cluster
extremely close to the unit circle as $N$ increases. In addition, the coefficients of $p(z)$ are all
small. (Note that the constant term $p_0$ of a monic Schur stable polynomial $p(z)=z^N+\cdots +p_0$
is always of modulus less than $1$.) For these two reasons, these polynomials are excellent test polynomials
for factorization algorithms.

\vsk
Finding nice Hurwitz stable polynomials is a much harder task. The Wilkinson polynomials
$p(z)=(z+1)(z+2) \cdots (z+N)$ have astronomically large coefficients and are therefore not
feasible. Well, one could take $p(z)=(z+1)^N$,
but already for $N=20$ the largest coefficient is $\tbinom{20}{10}=184756$. The choice
$p(z)=(z+\mu)^N$ with $0 < \mu < 1$ is also critical, since then the constant term $p_0=\mu^N$
may become the machine zero. Thus, I pose the following as a test: {\em find a Hurwitz stable
polynomial of degree 20 with positive integer coefficients such that the largest coefficient is about a hundred times better
than 184756, that is, such that it does not exceed 2000}.

\section{M\"obius transformation} \label{S2}

A point $z$ lies in the left open half-plane if and only if its distance to $-1$ is smaller than that to $1$,
that is, if and only if $|1-z|/|1+z| >1$. Consequently, if $u(z)$ has degree $N$ and all zeros of $u(z)$ are
of  modulus greater than $1$,
then $(1+z)^N u\left(\frac{1-z}{1+z}\right)$ is a Hurwitz stable polynomial of degree $N$.

\vsk
Let $u(z)=2+z+ \cdots+ z^N$ be the reverse of the polynomial we encountered in the introduction.
Then, for $z \neq 0$,
\begin{eqnarray}
\ell_N(z) & := & (1+z)^N u\left(\frac{1-z}{1+z}\right)=\frac{2(1+z)^{N+1}-(1-z)(1+z)^N-(1-z)^{N+1}}{2z}
\nonumber\\
& = & \frac{(1+z)^{N+1}-(1-z)^{N+1}}{2z}+(1+z)^N. \label{2.1}
\end{eqnarray}
Thus, $\ell_N(z)$ is a Hurwitz stable
polynomial of degree $N$. In what follows we frequently represent polynomials by their coefficient vectors as in Matlab,
that is, we write the polynomial $p_N z^N + \cdots +p_0$ as [$p_N$ ... $p_0$]. For even $N$, the coefficients of
the polynomials~(\ref{2.1}) are all even, and hence we divide them by $2$. The first polynomials are
\[{\tt [1 \; 3]}, \quad
{\tt [1  \;   1   \;  2]}, \quad
{\tt[1  \;   7  \;   3   \;  5]}, \quad
{\tt [1  \;   2   \;  8  \;   2  \;   3]}, \quad
{\tt [1\;  3\; 18\; 10\; 25\;  3\;  4]},\]
and for the degrees $10, 16, 20$ we obtain
\begin{eqnarray*}
N=10: & & {\tt [1\;  5\; 50\; 60\;270\;126\;336\; 60\;105\;  5\;  6]},\\
N=16: & & {\tt [1\; 8 \; 128 \; 280\;  2100\;  2184\; 10192 \; 5720 \; 18590 \; 5720 \; 13728\; 2184 \; 4004},\\
& &  {\tt \; 280 \; 400\; 8\; 9]}\\
N=20: & & {\tt [1\; 10\;200\;570\;  5415\; 7752\; 46512\; 38760\;164730\; 83980}\\
& & {\tt \; 268736\; 83980\;209950\; 38760\; 77520\;  7752\; 12597\; 570\;      760\; 10 \; 11]}.
\end{eqnarray*}
Figure 1 shows the zeros of $\ell_{20}(z)$.

\begin{figure}[thb] \label{Mob}
  \begin{center}
    \includegraphics[width=12cm]{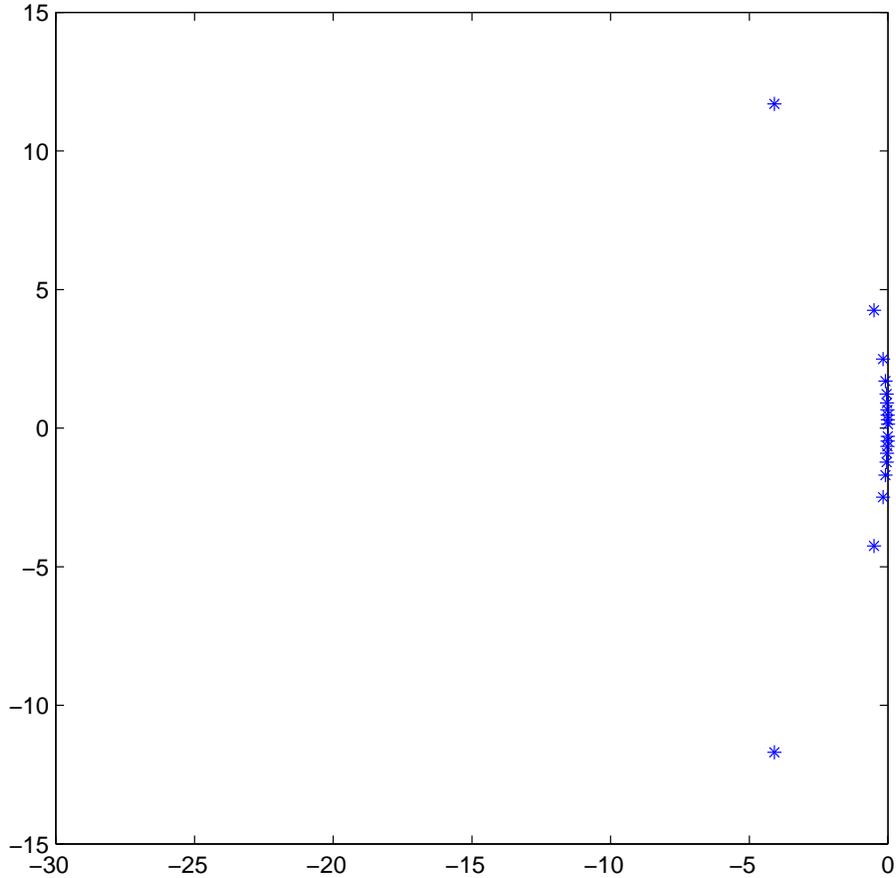}
    \caption{Zeros of polynomial (\ref{2.1}) for $N=20$.}
\end{center}
\end{figure}

\section{The inequality of Beauzamy} \label{S3}

The following is a slight improvement of an inequality which, for $v=1$, is stated (without proof) and attributed to
Beauzamy in~\cite{Trev}.

\begin{thm} \label{Theo 3.1}
Let $p_N(z)=p_N z^N+\cdots+p_0$ be a Hurwitz stable polynomial of even degree with $p_N \ge 1$ and $p_0 \ge 1$.
Then for every real number $v \ge 1$,
\[p_N(v)\ge (v^2+1)^{N/2}.\]
\end{thm}

{\em Proof.} Suppose $p_N(z)$
has exactly $n$ real zeros
and exactly $m$ pairs of genuinely complex conjugate zeros, multiplicities taken into account.
Then $N=n+2m$ and we may write
\begin{equation}
p_N(z)=p_N\prod_{j=1}^n(z+\mu_j)\prod_{j=1}^m (z^2+2 x_j z+  |w_j|^2) \label{3.1}
\end{equation}
with
$\mu_1, \ldots, \mu_n >0$ and with $w_j=-x_j+\ii y_j$, $x_j >0$ ($j=1, \ldots, m$).
The constant term of~(\ref{3.1}) is $p_0=p_N\mu_1 \cdots \mu_n |w_1|^2 \cdots |w_m|^2$.
Clearly, all coefficients of~(\ref{3.1})
are positive.
We have
\begin{eqnarray*}
p_N(v) & = & p_N\prod_{j=1}^n(v+\mu_j)\prod_{j=1}^m (v^2+2 x_j v+  |w_j|^2)\\
& \ge & \frac{p_N}{v^n}\prod_{j=1}^n(v^2+\mu_j v)\prod_{j=1}^m (v^2+ |w_j|^2).
\end{eqnarray*}
Put $\ph_1=\mu_1 v, \ldots, \ph_n=\mu_n v, \ph_{n+1}=|w_1|^2, \ldots, \ph_{n+m}=|w_m|^2$
and let $M=n+m$. Then
\[p_N(v) \ge \frac{p_N}{v^n}\prod_{j=1}^M(v^2+\ph_j), \quad \prod_{j=1}^M \ph_j =\frac{v^n}{p_N} p_0.\]
It follows that
\[p_N(v) \ge\frac{p_N}{v^n} \Big((v^2)^M+(v^2)^{M-1}s_1+(v^2)^{M-2}s_2+\cdots+(v^2)^{M-M}s_M\Big),\]
where
\begin{eqnarray*}
& & s_1 = \ph_1+\ph_2 + \cdots + \ph_M,\\
& & s_2 = \ph_1\ph_2 + \ph_1 \ph_3 + \cdots + \ph_{M-1}\ph_M,\\
& & \ldots, \\
& & s_M = \ph_1 \ph_2 \ldots \ph_M
\end{eqnarray*}
are the symmetric functions of $\ph_1, \ldots, \ph_M$.
 The sum $s_k$ ($1 \le k \le M-1$) contains $\tbinom{M}{k}$ terms
and each $\ph_j$ occurs exactly $\tbinom{M-1}{k-1}$ times in $s_k$.
The inequality between the arithmetic and geometric means therefore gives
\begin{eqnarray*}
s_k & \ge  & \binom{M}{k} \left(\ph_1^{M-1 \choose k-1} \ph_2^{M-1 \choose k-1}
\cdots \ph_M^{M-1 \choose k-1}\right)^{1\big/{M \choose k}}\\
& = & \binom{M}{k}(\ph_1\ph_2 \ldots \ph_M)^{k/M}=\binom{M}{k}s_M^{k/M}.
\end{eqnarray*}
Thus, using the binomial theorem and taking into account that $p_N \ge 1$, $p_0 \ge 1$, $v \ge 1$
we get
\begin{eqnarray*}
p_N(v) & \ge & \frac{p_N}{v^n} \sum_{k=0}^M \binom{M}{k}(v^2)^{M-k}s_M^{k/M}=\frac{p_N}{v^n}\left(v^2+s_M^{1/M}\right)^M\\
& = & \frac{p_N}{v^n}\left(v^2+\left(\frac{v^n p_0}{p_N}\right)^{1/M}\right)^M=\frac{1}{v^n}\left(p_N^{1/M}v^2+v^{n/M}p_0^{1/M}\right)^M\\
& \ge & \frac{1}{v^n}\left(v^2+v^{n/M}\right)^M \ge \frac{1}{v^n}(v^2+1)^M=\left(\frac{v^2+1}{v^2}\right)^{n/2}(v^2+1)^{n/2+m}\\
& \ge & (v^2+1)^{n/2+m} = (v^2+1)^{N/2}. \quad \square
\end{eqnarray*}

\begin{cor}\label{Cor 3.2}
Let $p_N(z)=p_N z^N+\cdots+p_0$ be a Hurwitz stable polynomial of even degree with $p_N \ge 1$ and $p_0 \ge 1$.
Then the sum of the coefficients is greater than or equal to $2^{N/2}$ and
at least one of the coefficients is greater than or equal to $2^{N/2}/(N+1)$.
\end{cor}

{\em Proof.} The sum of the coefficients is $p_N(1)$, and this is at least $2^{N/2}$ by Theorem~\ref{Theo 3.1}
with $v=1$. The polynomial has $N+1$ coefficients,
and denoting the maximum of the coefficients by $p_{\max}$, we have $p_N(1)\le (N+1)p_{\max}$,
which implies the asserted estimate $p_{\max} \ge 2^{N/2}/(N+1)$. $\;\:\square$

\vsg
The previous corollary provides us with a very crude lower bound for the largest coefficient. I conjecture that the $N+1$ can be replaced
by its square root, possibly with a multiplicative constant. However, this is not the point for our
purpose.

\begin{exa} \label{Exa 3.3}
{\rm
Suppose $p_{50}(z)=p_{50}z^{50}+\cdots+p_0$ is a Hurwitz stable polynomial of degree $50$.
Since a polynomial $p_N(z)$ of degree $N$ is Hurwitz stable if and only if so is the reverse polynomial $z^N p(1/z)$
(a property which is not shared by Schur stability), we may without loss of generality assume
that $p_{50} \le p_0$. Then we may write
\[p_{50}(z)=p_{50}\left(z^{50}+\frac{p_{49}}{p_{50}}z^{49}+ \cdots+\frac{p_{0}}{p_{50}}\right)\]
and apply Corollary \ref{Cor 3.2} to the polynomial in parentheses. We have $2^{25}=33554432$. Consequently, if $p_{\max}$
is the largest coefficient, then $p_{\max}/p_{50} \ge 2^{25}/51>650000$. If the coefficients are
required to be integers, this means that $p_{\max}> 650000$.

\vsk
The first even $N$ for which $2^{N/2}/(N+1) > 10000$ is $N=38$. Thus, a Hurwitz stable polynomial
of even degree
with positive integer coefficients not exceeding $10000$ must have a degree of at most
$36$. $\;\:\square$
}
\end{exa}

\begin{exa} \label{Exa 3.4}
{\rm
Let $p_{2}(z)=z^2+2 x z+ 1$.
Take $x=0.1$ and consider the polynomial $p_{20}(z):=p_{2}(z)^{10}=(z^2+2 x z+1)^{10}$.
The sum of the coefficients of $p_{20}(z)$ is about $2656$, which is comparable to $2^{N/2}=2^{10}=1024$.

\vsk
The polynomial $q_{20}(z)$ resulting from $p_{20}(z)$
by taking only the first $4$ digits of the coefficients after the comma is
\begin{eqnarray*}
& & {\tt [1.0000\;\;2.0000\;\; 11.8000\;\; 18.9600\;\; 59.7360\;\; 78.8006 \;\;  172.4294 \;\;  188.5647}\\
& & {\tt \;  315.8939 \;\;  286.4110 \;\;  384.8009 \;\;  286.4110 \;\;  315.8939 \;\;  188.5647 \;\;   172.4294}\\
& & {\tt \; 78.8006\;\; 59.7360 \;\; 18.9600\;\; 11.8000\;\; 2.0000\;\; 1.0000]}.
\end{eqnarray*}
Note that $q_{20}(z)$ has moderately sized coefficients which, in contrast to those of
$p_{20}(z)$, are precisely given within the machine precision. Figure 2 shows the zeros of $p_{20}(z)$
and $q_{20}(z)$ as they are given by Matlab.
Thus, the polynomial $q_{20}(z)$ has six zeros in the right half-plane and is therefore not Hurwitz stable!
$\;\: \square$
}
\end{exa}

\begin{figure}[thb] \label{Rund}
  \begin{center}
    \includegraphics[width=12cm, height=12cm]{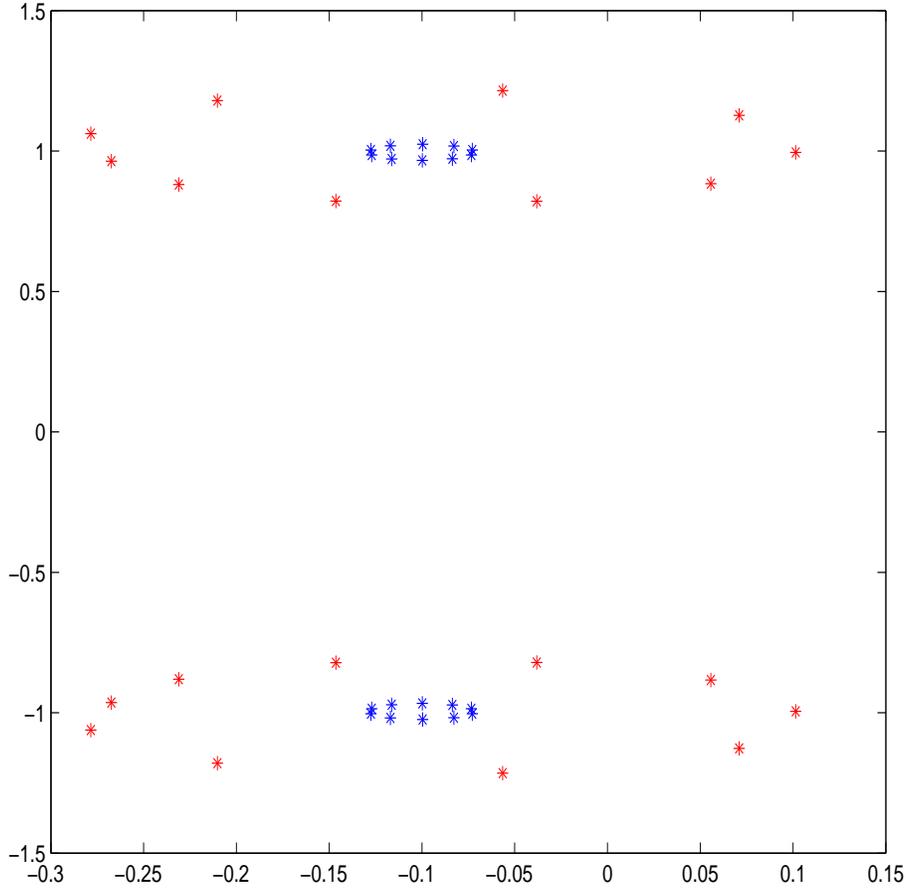}
   \caption{Zeros of $p_{20}(z)$ (blue) and of $q_{20}(z)$ (red) from Example \ref{Exa 3.4}
    obtained by Matlab.}
\end{center}
\end{figure}

\section{Integer coefficients} \label{S4}

We consider Hurwitz stable polynomials $p(z)$ whose coefficients are positive integers.
The degree is denoted by $N$, the maximal coefficient by $p_{\max}$, the sum of
the coefficients by $\si$, and the maximum of the real parts of the roots ($=$ spectral abscissa) by $\al$.
For each degree $N$, there
are two kinds of optimal polynomials: the polynomials with minimal largest coefficient and
the polynomials with minimal sum $\si$ of the coefficients.  We call these polynomials $c$-optimal and $\si$-optimal.
Small degrees $N$ are easy, because all possible cases can be checked by Matlab.

\vsg
$\boldsymbol{ N = 1.}$ The polynomial $a_{1}(z)=z+1={\tt [1 \;1]}$ is the best in all respects.

\vsg
$\boldsymbol{ N = 2.}$ Here the polynomial $a_{2}(z)=z^2+z+1={\tt [1 \; 1 \; 1]}$ is optimal on all accounts.
Its spectral abscissa is $-1/2$.

\vsg
$\boldsymbol{ N = 3.}$ The five $c$-optimal polynomials
are
\begin{eqnarray*}
& & b_{3}(z)={\tt [1 \; 1\;  2\;  1]}, \quad \al= -0.2151, \quad \si=  5,\\
& & c_{3}(z)={\tt [1\;  2\;  1\;  1]}, \quad \al=    -0.1226, \quad \si=   5,\\
& & d_{3}(z)={\tt [1\;  2\;  2\;  1]}, \quad \al=    -0.5000, \quad \si=   6,\\
& & e_{3}(z)={\tt [1\;  2\;  2\;  2]}, \quad \al=    -0.2282, \quad \si=   7,\\
& & f_{3}(z)={\tt [2\;  2\;  2\;  1]}, \quad \al=     -0.1761, \quad \si=   7,
\end{eqnarray*}
and the first two of them are the $\si$-optimal polynomials.

\vsg
$\boldsymbol{ N = 4.}$ The list of the nine $c$-optimal polynomials is
\begin{eqnarray*}
& & a_4(z)={\tt [1\;  1\;  3\;  1\;  1]}, \quad \al=  -0.1484, \quad \si=    7  ,\\
& & c_4(z)={\tt [1\;  1\;  3\;  2\;  1]}, \quad \al=   -0.1049, \quad \si=    8 ,\\
& & d_4(z)={\tt [1\;  2\;  3\;  1\;  1]}, \quad \al=    -0.0433, \quad \si=    8,\\
& & e_4(z)={\tt [1\;  2\;  3\;  2\;  1]}, \quad \al=    -0.5000, \quad \si=    9,\\
& & f_4(z)={\tt [1\;  2\;  3\;  3\;  1]}, \quad \al=    -0.2151, \quad \si=   10,\\
& & g_4(z)={\tt [1\;  3\;  3\;  2\;  1]}, \quad \al=   -0.1226, \quad \si=   10 ,\\
& & h_4(z)={\tt [1\;  2\;  3\;  3\;  2]}, \quad \al=   -0.0433, \quad \si=   11 ,\\
& & i_4(z)={\tt [1\;  3\;  3\;  3\;  1]}, \quad \al=    -0.1910, \quad \si=   11,\\
& & j_4(z)={\tt [2\;  3\;  3\;  2\;  1]}, \quad \al=   -0.0287, \quad \si=   11 .
\end{eqnarray*}
The only $\si$-optimal polynomial is $a_4(z)$.

\vsg
$\boldsymbol{ N = 5.}$ The nine $c$-optimal polynomials are
\begin{eqnarray*}
& & b_5(z)={\tt [1 \; 1\;  4\;  3\;  2\;  1]}, \quad \al=  -0.0835, \quad \si=  12,\\
& & c_5(z)={\tt [1\;  2\;  3\;  4\;  1\;  1]}, \quad \al=    -0.0320, \quad \si=   12,\\
& & d_5(z)={\tt [1\;  1\;  4\;  3\;  3\;  1]}, \quad \al=    -0.0582, \quad \si=   13,\\
& & e_5(z)={\tt [1\;  2\;  3\;  4\;  2\;  1]}, \quad \al=    -0.0354, \quad \si=   13,\\
& & f_5(z)={\tt [1\;  2\;  4\;  3\;  2\;  1]}, \quad \al=     -0.0204, \quad \si=   13,\\
& & g_5(z)={\tt [1\;  3\;  3\;  4\;  1\;  1]}, \quad \al=   -0.0203, \quad \si=   13,\\
& & h_5(z)={\tt [1\;  2\;  4\;  4\;  2\;  1]}, \quad \al=     -0.1484, \quad \si=   14,\\
& & i_5(z)={\tt [1\;  2\;  4\;  4\;  3\;  1]}, \quad \al=     -0.2151, \quad \si=   15,\\
& & j_5(z)={\tt [1\;  3\;  4\;  4\;  2\;  1]}, \quad \al=    -0.1226, \quad \si=   15,
\end{eqnarray*}
and the $\si$-optimal polynomials are $b_5(z)$ and $c_5(z)$. Of course, it might be that there
exist $\si$-optimal polynomials with $p_{\max} \ge 5$ and $\si \le 11$. However, the
coefficients of such polynomials are either a permutation of $6,1,1,1,1,1$ or a permutation
of $5, 2, 1,1,1,1$, and none of these $36$ polynomials is Hurwitz stable.

\vsg
$\boldsymbol{ N = 6.}$ We have the five $c$-optimal polynomials
\begin{eqnarray*}
& & b_6(z)={\tt [1\;  1\;  5\;  3\;  5\;  1\;  1]}, \quad \al= -0.0485, \quad \si=   17,\\
& & c_6(z)={\tt [1\;  1\;  5\;  4\;  5\;  2\;  1]}, \quad \al= -0.0393, \quad \si=    19 ,\\
& & d_6(z)={\tt [1\;  2\;  4\;  5\;  4\;  2\;  1]}, \quad \al=  -0.0399, \quad \si=    19,\\
& & e_6(z)={\tt [1\;  2\;  5\;  4\;  5\;  1\;  1]}, \quad \al=  -0.0108, \quad \si=    19 ,\\
& & f_6(z)={\tt [1\;  2\;  5\;  5\;  5\;  2\;  1]}, \quad \al=  -0.1484, \quad \si=    21 ,
\end{eqnarray*}
and the first of them is the only $\si$-optimal polynomial.

\vsg
$\boldsymbol{ N = 7.}$ Beginning with this degree things become challenging. Inspection of
the $7^8=5764801$ polynomials with $p_{\max} \le 7$ shows that exactly two of them are Hurwitz stable,
namely,
\begin{eqnarray*}
& & b_7(z)={\tt [1\;  2\;  5\;  7\;  7\;  6\;  2\;  1]}, \quad \al= -0.0175, \quad \si=   31,\\
& & c_7(z)={\tt [1\;  2\;  6\;  7\;  7\;  5\;  2\;  1]}, \quad \al= -0.0077, \quad \si=    31 .
\end{eqnarray*}
Consequently, these are the $c$-optimal polynomials of degree $7$. Note that each of the two polynomials is
the reverse of the other one. These two polynomials are not $\si$-optimal, because, for example, we also have the polynomials
\begin{eqnarray*}
& & d_7(z)={\tt [1 \;  2 \;  5 \;  8 \;  5 \;  6 \;  1 \;  1]}, \quad \al=-0.0131, \quad \si=29,\\
& & e_7(z)={\tt [1 \;  3 \;  4 \;  9\;  4 \;  6 \;  1 \;  1]}, \quad \al=-0.0526, \quad \si=29.
\end{eqnarray*}
I have not examined whether the last two polynomials are $\si$-optimal.

\vsg
{\bf Multiplication and doubling.} One way of getting Hurwitz stable polynomial of higher degrees is to multiply Hurwitz stable polynomials
of lower degrees. Another way is as follows. Since $z$ is in the open left half-plane if and only if so is $z+1/z$, it follows
that $q_N(z)$ is Hurwitz stable of degree $N$ if and only if $p_{2N}(z)=z^Nq_N(z+1/z)$ is a Hurwitz stable polynomial of degree $2N$.
We refer to the passage from $q_N(z)$ to $p_{2N}(z)=z^Nq_N(z+1/z)$ as doubling. The sum of the coefficients
of a product is equal to the product of the sums of the coefficients, and if $p_{2N}(z)$ results from $q_N(z)$ by doubling,
then the sum of the coefficients of $p_{2N}(z)$ is $p_{2N}(1)=q_N(2)$.

\vsg
$\boldsymbol{ N = 8.}$  Multiplying $a_{4}(z)$ by itself we obtain
\[b_8(z)={\tt [1\;  2\;  7\;  8\; 13\; 8\;  7\;  2\;  1]}, \quad \al=-0.1484, \quad \si=7\cdot 7=49.\]
However, such products are usually far away from the optimal polynomials.
Doubling of $a_4(z)$ gives
\[a_8(z)={\tt [1 \; 1 \; 7 \; 4 \; 13 \; 4 \; 7 \; 1 \; 1]}, \quad \al=-0.0518, \quad \si= 39.\]
I don't know whether $a_8(z)$ is $\si$-optimal. It is surely not $c$-optimal, because
the polynomials
\begin{eqnarray*}
& & c_8(z)={\tt [1\;  2\;  6\;  9\; 11\; 10\;  7\;  2\;  1]},\quad \al= -0.0171,\quad \si=49,\\
& & d_8(z)={\tt [1\;  2\;  6\;  9\; 11\; 11\;  7\;  3\;  1]},\quad \al= -0.0075, \quad \si= 51,\\
& & e_8(z)={\tt [1\;  2\;  7\; 11\; 11\; 11\;  6\;  3\;  1]},\quad \al= -0.0135, \quad \si=53
\end{eqnarray*}
are Hurwitz stable. Note that the reverses of these polynomials are Hurwitz stable, too.
Clearly, these six polynomials are closer to the $c$-optimal polynomials.

\vsg
$\boldsymbol{ 10 \le N \le 18.}$ Doubling of $c_5(z)$, $d_6(z)$, $d_7(z)$, $c_8(z)$ yields
%& & N=10, \quad b_{10}(z)= {\tt [1 \;  1 \;  9 \;  7 \; 24 \; 13 \; 24 \;  7 \;  9 \;  1 \;  1]},\quad \al=-0.0172, \quad \si=97,\\
\begin{eqnarray*}
& & N=10, \quad b_{10}(z)={\tt [1\;  2\;  8\; 12\; 20\; 21\; 20\; 12\;  8\;  2\;  1] }, \quad \al= -0.0117, \quad \si= 107 ,\\
& & N=12, \quad b_{12}(z)={\tt [1\;  2\; 10\; 15\; 35\; 37\; 53\; 37\; 35\; 15\; 10\;  2\;  1] },\\
& &  \qquad \qquad \quad \al= -0.0134 , \quad \si= 253 ,\\
& & N=14, \quad b_{14}(z)={\tt [1\;  2\; 12\; 20\; 51\; 68\;101\;101\;101\; 68\; 51\; 20\; 12\;  2\;  1] },\\
& & \qquad \qquad \quad \al= -0.0050, \quad \si= 611 ,\\
& & N=16, \quad b_{16}(z)={\tt [1\;  2\; 14\; 23\; 75\; 97\;197\;192\;271\;192\;197\; 97\; 75\; 23\; 14\;  2\;  1] },\\
& & \qquad \qquad \quad \al= -0.0042, \quad \si= 1473 ,
\end{eqnarray*}
and for $N=18$, I found
\[b_{18}(z)={\tt [1\;  2\; 16\; 27\; 98\;139\;303\;353\;523\;479\;523\;353\;303\;139\; 98\; 27\; 16\;  2\;  1] }\]
with $\al= -0.0046$ and $\si= 3403$.

\vsg
$\boldsymbol{ N = 20.}$ We arrive at the test posed in the introduction. Multiplication of $b_{10}(z)$ by itself
gives
\[{\tt [1\; 4\;    20\;    56\;   152\;   314\;   588\;   920\;  1288\;  1548\;  1667\;  1548\;  1288
\;   920\;   588\;   314\;   152\;    56\;    20\;     4\;     1]}\]
with $\al= -0.0117$ and $\si=11449$. Multiplying other combinations yields similar results, the best being
$a_{2}(z)\cdot b_{18}(z)$, which is
\[{\tt [1\;     3\;    19\;    45\;   141\;   264\;   540\;   795\;  1179\;  1355\;  1525\;  1355\;  1179
\;   795\;   540\;   264\;   141\;    45\;    19\;     3\;     1]}\]
with $\al=-0.0046$ and $\si=10209$.
Thus, eventually we easily passed the test and constructed
a polynomial with  $p_{\max} = 1525$. However, notice that the success resulted from knowing the very good
polynomials $b_{10}(z)$ and $b_{18}(z)$. In fact we can do it even better.
Doubling $b_5(z)={\tt [1 \;   1  \;  4 \; 3  \;  2 \; 1]}$ we get
\[c_{10}(z)={\tt [1 \; 1 \; 9 \; 7 \;   24  \;  13 \;   24 \; 7 \; 9 \; 1 \; 1]}, \quad \si=97,\]
and doubling this again, we arrive at
\begin{equation}
{\tt [1 \; 1 \; 19 \; 16 \; 141 \; 98 \; 540  \;  303 \;   1179 \; 523 \; 1525 \; 523  \;  1179 \;
303 \;   540 \;   98  \;  141 \; 16 \; 19 \; 1 \; 1]},\label{c20}
\end{equation}
with $\al=-0.0067$ and $\si=7167$,
which has the smallest $\si$ we have found. This polynomial will be denoted by $c_{20}(z)$.
Moreover, doubling of $b_{10}(z)$ yields
the polynomial
\[
{\tt [1\; 2 \;18\; 30\; 129\;177\;484\; 537\;1046\;920\;1349\;920\;1046\;537\;484\;177\; 129\; 30\; 18\;2 \;1]}\]
with $p_{\max}=1349$, $\al=-0.0038$, $\si=8037$. We henceforth denote this polynomial by $b_{20}(z)$.
I have not found a Hurwitz stable polynomial of degree $20$ whose largest coefficient is smaller than $1349$.
Figures~3 and~4 show the zeros of $b_{18}(z)$ and $b_{20}(z)$.

\begin{figure}[thb] \label{Mult}
  \begin{center}
    \includegraphics[width=12cm]{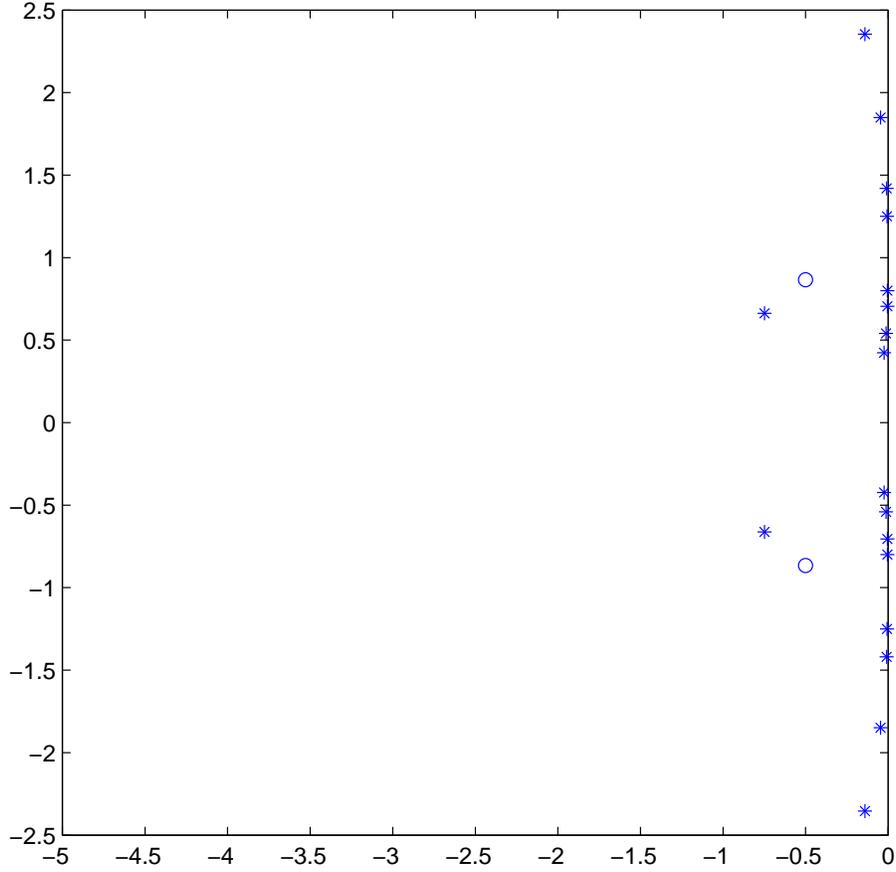}
    \caption{Location of the zeros of the product of the polynomials $a_{2}(z)$ (two circles)
    and $b_{18}(z)$ (18 asterisks).}
\end{center}
\end{figure}

\begin{figure}[thb] \label{Rec}
  \begin{center}
    \includegraphics[width=12cm]{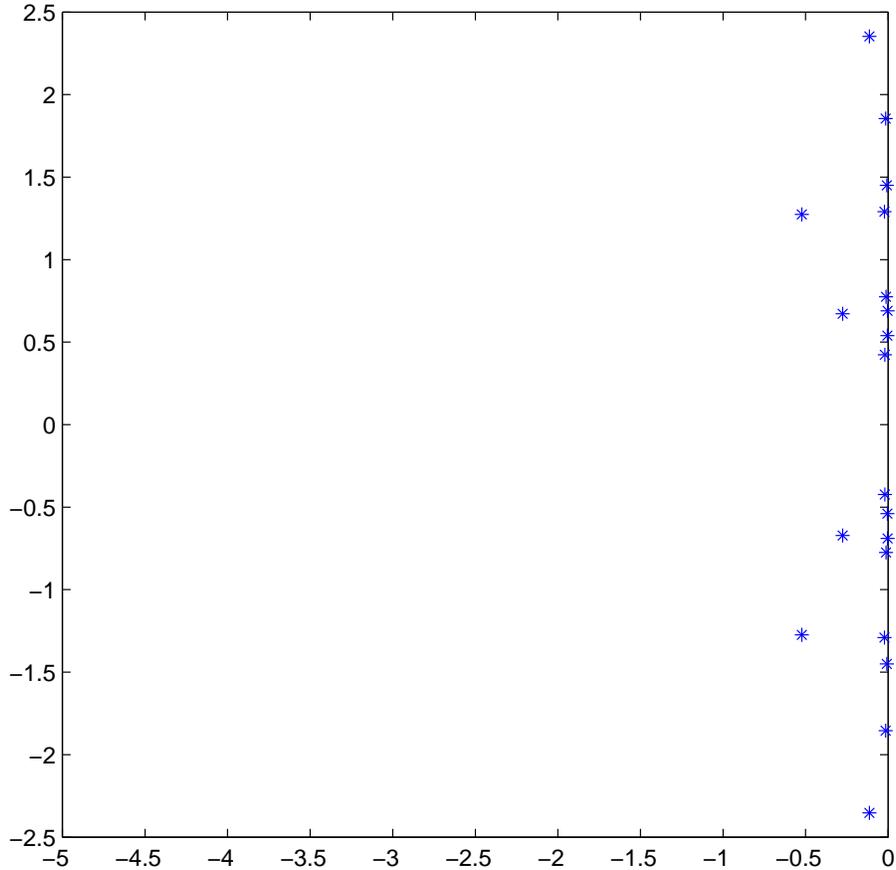}
    \caption{Zeros of the polynomial $b_{20}(z)$.}
\end{center}
\end{figure}

{\bf Powers of $\boldsymbol{2}$.} Repeated doubling of $a_1(z)={\tt [1\; 1]}$ yields the polynomials
\begin{eqnarray*}
& & a_{2}(z)={\tt [1 \; 1 \; 1]}, \quad \si=3,\\
& & a_4(z)={\tt [1\;  1\;  3\;  1\;  1]}, \quad \si=    7,\\
& & a_8(z)={\tt [1 \; 1 \; 7 \; 4 \; 13 \; 4 \; 7 \; 1 \; 1]}, \quad \si= 39,\\
& & a_{16}(z)=[1\;1\;15\;11\;83\;45\;220\;88\;303\;88\;220\;45\;83\;11\; 15\;1\;1], \quad \si=1231,
\end{eqnarray*}
and finally
\begin{eqnarray*}
& & a_{32}(z)={\tt [1 \; 1 \; 31 \; 26 \; 413 \; 293 \; 3141 \; 1896 \; 15261 \; 7866 \; 50187 \; 22122 \; 115410 \; 43488}\\
& & \qquad \qquad \;{\tt  189036 \; 60753 \; 222621 \;
60753 \; 189036 \; 43488 \; 115410 \; 22122 \; 50187 \; 7866}\\
& & \qquad \qquad  \; 1{\tt 5261 \; 1896 \; 3141 \; 293 \; 413 \; 26 \; 31 \; 1 \; 1]}, \quad \si=1242471.
\end{eqnarray*}
The sum of the coefficients of $a_{64}(z)$ equals $a_{32}(2) \approx 1.2791\cdot 10^{12}$. When comparing  $a_{16}(z)$
and $b_{16}(z)$, we see that $a_{16}(z)$ has the smaller $\si$ and that $b_{16}(z)$ has the smaller $p_{\max}$.

\section{Higher degrees} \label{S5}

Let $N$ be even. From Theorem \ref{Theo 3.1} we infer that if $p_N(z)=p_Nz^N+\cdots+p_0$ is Hurwitz stable
with $p_N \ge 1$ and $p_0 \ge 1$, then the sum of the coefficients always satisfies $p_N(1) \ge 2^{N/2}$.
The polynomials $p_N(z)=(z^2+2xz+1)^{N/2}$ with sufficiently small $x >0$ show that Theorem~\ref{Theo 3.1} is sharp.
In particular, given any $\eps >0$, there is such a polynomial for which $p_N(1)< (2+\eps)^{N/2}$.
But what happens if the coefficients are required to be integers?

\vsk
Let $p_{\max}(N)$ denote the minimum of
the largest coefficients and let $\si(N)$ be the minimum of the sum of the coefficients of the Hurwitz stable polynomials
of degree $N$ with positive integer coefficients. Equivalently, $p_{\max}(N)$ is the largest coefficient of the $c$-optimal
polynomials and $\si(N)$ is the sum of the coefficients of the $\si$-optimal polynomials.
From the previous section we know the following.

\[\begin{array}{|c||c|c|c|c|c|c|c|c|c|c|c|}
\hline
N & 1 & 2 & 3 & 4 & 5 & 6 & 7 & 8  & 10 & 20 & 32\\
\hline
p_{\max}(N) & 1 & 1 & 2 & 3 & 4 & 5 & 7 & 11 & \le 21 & \le 1349 & \le 222621\\
\hline
\si(N) & 2 & 3 & 5 & 7 & 12 & 17 & \le 29 & \le 39  & \le 97 & \le 7167 & \le 1242471\\
\hline
\end{array}\]

\begin{thm} \label{Theo Fek}
The limits of $\si(N)^{1/N}$ and $p_{\max}(N)^{1/N}$ exist and we have
\[\lim_{N\to \iy}\si(N)^{1/N}=\inf_{N \ge 1}\si(N)^{1/N}=\lim_{N\to \iy}p_{\max}(N)^{1/N}
=\inf_{N \ge 1}p_{\max}(N)^{1/N}=:\be.\]
\end{thm}

{\em Proof.}
Since the product of Hurwitz stable polynomials is again Hurwitz stable, it follows that
\[\si(N+M) \le \si(N)\si(M).\]
Fekete's subadditive lemma, for which see \cite{Fek}\footnote[1]{Interestingly, in their equally titled
papers \cite{Fek}, \cite{Schur},
Schur and Fekete considered the problem whether there are infinitely many polynomials with
integer coefficients and given leading coefficient whose zeros are all simple and lie in a compact subset $E$
of the plane. For example, in the case where $E$ is a half-disk with diameter $2R<3\sqrt{3}/2=2.5981\dots$,
Fekete showed that the number of such polynomials must be finite. We here are concerned
with the case where $E$ is the open left half-plane, which is neither bounded nor closed.
And indeed, the polynomials $(z+j)\ldots(z+j+N)$ ($j=1,2,\ldots$) constitute an infinite family of monic polynomials
of even fixed degree
with integer coefficients whose zeros are all simple and are located in $E$.}
or \cite[p.~16]{Law}, therefore implies that the limit of $\si(N)^{1/N}$
exists and coincides with the infimum of $\si(N)^{1/N}$ for $N \ge 1$. For every polynomial $p_N(z)$
of degree $N$ with positive coefficients, the inequalities
\[\frac{\si}{N+1} \le p_{\max} \le \si\]
hold, where $\si$ is the sum and $p_{\max}$ is the maximum of the coefficients. If $p_N(z)$ is $\si$-optimal,
then $\si=\si(N)$ and hence
\[p_{\max}(N) \le p_{\max} \le \si =\si(N).\]
In case  $p_N(z)$ is $c$-optimal, we have $p_{\max}=p_{\max}(N)$ and consequently,
\[\frac{\si(N)}{N+1} \le \frac{\si}{N+1} \le p_{\max}=p_{\max}(N).\]
Thus, $\si(N)/(N+1) \le p_{\max}(N) \le \si(N)$, which shows that the limit and the infimum of $p_{\max}(N)^{1/N}$
coincide with the limit and the infimum of $\si(N)^{1/N}$: $\;\:\square$

\vsg
Theorem~\ref{Theo 3.1}
with $v=1$ shows that $\be \ge \sqrt{2}=1.4142\ldots$.

\begin{prop} \label{Prop Wu}
Let $p_k(z)$ be any Hurwitz stable polynomial with positive integer coefficients. If $N$ is divisible by $k$,
then $\si(N) \le (\sqrt[k]{p_k(1)}\,)^N$.
\end{prop}

{\em Proof.} Let $N=nk$ and consider $p_N(z)=p_k(z)^n$. Then
\[\si(N) \le p_N(1) = p_k(1)^n = (p_k(1)^{1/k})^N. \quad \square \]

The best results from Proposition \ref{Prop Wu} are delivered by taking $\si$-optimal polynomials,
in which case $p_k(1)=\si(k)$. Here are the numbers.

\[\begin{array}{|c||c|c|c|c|c|c|c|c|c|c|c|}
\hline
k & 1 & 2 & 3 & 4 & 5 & 6 & 7 & 8  & 10 & 20 & 32\\
\hline
\si(k) \le & 2 & 3 & 5 & 7 & 12 & 17 & 29 & 39  & 97 & 7167 & 1242471\\
\hline
 &  &  &  &  &  &  &  &   &  &  & \\[-2.4ex]
\sqrt[k]{\si(k)}\le  & 2 & 1.74 & 1.72 & 1.63 & 1.65 &   1.61 &   1.62 & 1.59 & 1.59 & 1.56 & 1.56\\
\hline
\end{array}\]

Since $\sqrt[20]{7167}=1.5587\ldots $ and $\sqrt[32]{1242471}=1.5504 \ldots$, we arrive at the following.

\begin{cor} \label{Cor be1}
We have
\[1.4142\ldots = \sqrt{2} \le \be \le \sqrt[32]{1242471}=1.5504\ldots.\]
\end{cor}

\begin{cor} \label{Theo 5.1}
If $N$ is divisible by $20$ or $32$, then $1.41^N < \si(N) < 1.56^N$.
\end{cor}

{\em Proof.} We know from Theorem \ref{Theo 3.1} that $\si(N) \ge  2^{N/2}=(2^{1/2})^N > 1.41^N$,
and Proposition~\ref{Prop Wu} implies that $\si(N) \le (\sqrt[20]{7167}\,)^N < 1.56^N$ if $N$ is divisible by $20$
and that $\si(N) \le (\sqrt[32]{1242471}\,)^N < 1.56^N$ if $N$ is divisible by $32$. $\;\:\square$

\vsg
In what follows we need the sequence $v_0,v_1, v_2, \ldots$ given by $v_0=1$ and $v_{n+1}=v_n+1/v_n$.
The first terms are
\[v_0=1, \quad v_1=2, \quad v_2= \frac{5}{2}=2.5, \quad v_3= \frac{29}{10}=2.9, \quad
v_4= \frac{941}{290}=3.2448\ldots .\]

\begin{lem}\label{Lem vn}
We have $\sqrt{n+1} < v_n < 2\sqrt{n}$ for $n \ge 2$ and
\[v_n=\sqrt{2n}\left(1+O\left(\frac{\log n}{n}\right)\right).\]
\end{lem}

{\em Proof.}
We prove the inequalities $\sqrt{n+1} < v_n \le 2\sqrt{n}$ by induction on $n$. They are
obviously true for $n=2$. So suppose they hold for $n$. We then have
\[v_{n+1}=v_n+\frac{1}{v_n} \le 2\sqrt{n}+\frac{1}{\sqrt{n+1}} < 2\sqrt{n+1}\]
because
\[2\sqrt{n+1}-2\sqrt{n}=\frac{2}{\sqrt{n+1}+\sqrt{n}} > \frac{2}{\sqrt{n+1}+\sqrt{n+1}}=\frac{1}{\sqrt{n+1}}.\]
In the same vein,
\[v_{n+1}=v_n+\frac{1}{v_n} > \sqrt{n+1}+\frac{1}{2\sqrt{n}} > \sqrt{n+2}\]
since
\[\sqrt{n+2}-\sqrt{n+1} = \frac{1}{\sqrt{n+2}+\sqrt{n+1}}< \frac{1}{\sqrt{n}+\sqrt{n}}=\frac{1}{2\sqrt{n}}.\]
This completes the proof of the inequalities.
To prove the asymptotics, note first that
the numbers $v_n$ satisfy $v_n^2=v_{n-1}^2+1/v_{n-1}^2+2$. Consequently,
\[v_n^2 > v_{n-1}^2+2 > v_{n-2}^2 +2\cdot 2 > \ldots > v_0^2 +n\cdot 2 =2n+1,\]
which implies that $v_n > \sqrt{2n}$. On the other hand,
\[v_n^2=v_{n-1}^2+\frac{1}{v_{n-1}^2}+2= v_{n-2}^2+\frac{1}{v_{n-2}^2}+\frac{1}{v_{n-1}^2}+2\cdot 2\]
and so on, which eventually gives
\[v_n^2=v_1^2+\frac{1}{v_1^2}+ \cdots+\frac{1}{v_{n-1}^2}+(n-1)\cdot 2.\]
As $v_k > \sqrt{2k}$, we conclude that
\[\frac{v_n^2}{2n}< \frac{1}{2n}\left(4+\frac{1}{2\cdot 1}+\cdots+\frac{1}{2(n-1)}\right)+\frac{n-1}{2n}\cdot 2
=1+O\left(\frac{\log n}{n}\right)\]
and hence
\[v_n\le\sqrt{2n}\left(1+O\left(\frac{\log n}{n}\right)\right).\]
This estimate in conjunction with the inequality $v_n > \sqrt{2n}$ proves the lemma. $\:\;\square$

\vsg
We define
\[\ga_k=\sum_{j=0}^{k-1}\frac{\log v_j}{2^{j+1}}, \quad \ga=\sum_{j=0}^\iy \frac{\log v_j}{2^{j+1}}=0.4329\ldots.\]
The first values are
\[\begin{array}{|c|c|c|}
\hline
k & \ee^{\ga_k} & \ga_k\\
\hline
1 & 1 & 0 \\
2 & \sqrt[4]{2}=1.1892\ldots & 0.1733\ldots\\
3 & \sqrt[8]{10}=1.3335\ldots & 0.2878\ldots \\
4 & \sqrt[16]{290}=1.4252\ldots & 0.3544 \ldots\\
\hline
\end{array}\]

Throughout what follows, if $N=2^n$ is a power of $2$, we denote by $a_N(z)$ the polynomials obtained from
$a_1(z)=z+1$ by $n$ doublings. The first of these polynomials are listed at the end of Section~\ref{S4}.

\begin{lem}\label{Lem aN1}
Let $N=2^n$. Then
\[a_N(1)=(v_n+1) \ee^{\ga_n N} < (2\sqrt{n}+1)\ee^{\ga N}\]
with $\ee^\ga= 1.5417\ldots$ and
\[a_N(1)=\left((2n)^{1/4}+(2n)^{-1/4}\right)\ee^{\ga N}\left(1+O\left(\frac{\log n}{n}\right)\right).\]
\end{lem}

{\em Proof.} We have
\[a_N(1)  =  v_0^{N/2}a_{n/2}(v_1)=v_0^{N/2}v_1^{N/4}a_{N/4}(v_2)
 =  \ldots = v_0^{N/2}v_1^{N/4} \cdots v_{n-1} a_1(v_n),\]
and since
\[\log( v_0^{N/2}v_1^{N/4} \cdots v_{n-1})  =  2^n\left(\frac{\log v_0}{2}+\frac{\log v_1}{2^2}+ \cdots + \frac{\log v_{n-1}}{2^n}\right)
=N\ga_n,\]
we get $a_N(1)=(v_n+1)\ee^{\ga_n N}$. The upper bound $(2\sqrt{n}+1)\ee^{\ga N}$ follows from
Lemma~\ref{Lem vn}. To prove the asymptotics, we write
\[\log a_N(1)=\log(v_n+1)+2^n\ga_n=\log(v_n+1)+2^n\ga-r_n, \quad r_n=\sum_{j=0}^\iy \frac{\log v_{n+j}}{2^{j+1}}.\]
By Lemma \ref{Lem vn},
\begin{equation}
r_n=\sum_{j=0}^\iy \frac{\log 2}{2^{j+2}}+\sum_{j=0}^\iy \frac{\log (n+j)}{2^{j+2}}
+O\left(\sum_{j=0}^\iy \frac{\log (n+j)}{2^{j+1}(n+j)}\right).\label{6.1}
\end{equation}
The first sum in (\ref{6.1}) is $(1/4)\log 2$. Using that $\log(1+x)<x$ for $x >0$, the second sum
can be estimated as follows:
\begin{eqnarray*}
\sum_{j=0}^\iy \frac{\log (n+j)}{2^{j+2}} & = & \sum_{j=0}^\iy \frac{\log n+\log(1+j/n)}{2^{j+2}}\\
& = & \frac{1}{4}\log n+O\left(\sum_{j=0}^\iy \frac{j}{2^{j+2}n}\right)=
\frac{1}{4}\log n+O\left(\frac{1}{n}\right).
\end{eqnarray*}
Multiplying the $j$th term in the third sum by $n/\log n$, it becomes
\[\frac{\log(n+j)}{\log n} \frac{n}{n+j} \frac{1}{2^{j+1}}
=\frac{\log n+\log(1+j/n)}{\log n} \frac{n}{n+j} \frac{1}{2^{j+1}}
< \left(1+\frac{j}{n \log n}\right)\frac{1}{2^{j+1}},\]
and as this is smaller than $(1+j)/2^{j+1}$, we arrive at the conclusion that
the third term in~(\ref{6.1}) is $O((\log n)/n)$. Putting things together we obtain
\[a_N(1)=(v_n+1)\ee^{\ga N}\ee^{-(1/4)\log n}\ee^{-(1/4)\log 2}\left(1+O\left(\frac{\log n}{n}\right)\right).\]
From Lemma \ref{Lem vn} we infer that
\[v_n+1=\sqrt{2n}\left(1+O\left(\frac{\log n}{n}\right)\right)+1
=(\sqrt{2n}+1)\left(1+O\left(\frac{\log n}{n}\right)\right).\]
What finally results is
\[a_N(1)=\frac{\sqrt{2n}+1}{(2n)^{1/4}}\ee^{\ga N}\left(1+O\left(\frac{\log n}{n}\right)\right),\]
which is equivalent to the assertion. $\;\:\square$

\vsg
Here is a slight improvement of Corollary \ref{Cor be1}.

\begin{cor} \label{Cor be2}
We have
\[1.4142\ldots = \sqrt{2} \le \be \le \ee^\ga=1.5417\ldots.\]
\end{cor}

{\em Proof.} Let $\eps >0$ be arbitrarily given. Choose $K=2^k$ so that
$(2\sqrt{k}+1)^{1/K} < 1+\eps$. Lemma~\ref{Lem aN1} then
gives $a_K(1)^{1/K} < (1+\eps )\ee^{\ga}$. If $N$ is divisible by $2^k$,
Proposition~\ref{Prop Wu} implies that $\si(N)^{1/N} \le a_K(1)^{1/K}< (1+\eps)\ee^{\ga}$,
whence $\be \le (1+\eps)\ee^\ga$. As $\eps >0$ was arbitrary, we conclude that
$\be \le \ee^\ga$. $\;\:\square$

\vsg
A polynomial $p_N(z)=p_Nz^N+\cdots+p_0$ of even degree $N$ is called symmetric if $p_j=p_{N-j}$ for all $j$. In that case
there is a unique polynomial $p_{N/2}(z)$ of degree $N/2$ such that $p_N(z)=z^{N/2} p_{N/2}(z+1/z)$.
The polynomial $p_{N/2}(z)$ is Hurwitz stable if and only if so is $p_N(z)$, and $p_{N/2}(z)$ has integer coefficients
if and only if $p_N(z)$ has integer coefficients.
If $N/2$ is also even and $p_{N/2}(z)$ is symmetric, we call $p_N(z)$ a $2$-fold symmetric polynomial.
We then have $p_{N/2}(z)=z^{N/4}p_{N/4}(z)$. If $N/4$ is even and $p_{N/4}(z)$ is symmetric, then $p_N(z)$
is said to be $3$-fold symmetric and so on. In other words, a polynomial is $k$-fold symmetric if and only
if it results after $k$ doubling procedures from another polynomial. Symmetry is $1$-fold symmetry in this context.

\vsk
We denote by $\si_k(N)$ the minimum of the sum of the coefficients among all Hurwitz stable $k$-fold symmetric polynomials
of degree $N$ with positive integer coefficients. Clearly,
$\si(N) \le \si_1(N) \le \si_2(N) \le \ldots$.

\begin{thm} \label{Theo kfold}
Let $N$ be divisible by $2^k$. Then
\[\si(N) \le \si_k(N) \le (v_k+1)^{N/2^k}\ee^{\ga_k N}\]
and
\[\ee^{\ga_k N} \le \Big((v_k^2+1)^{1/2^{k+1}}\Big)^N \ee^{\ga_k N}\le \si_k(N).\]
\end{thm}

{\em Proof.} Let $N=2^km$. The polynomial $p_N(z)=a_{2^k}(z)^m$ is $k$-fold symmetric and hence
$\si_k(N) \le p_N(1)$. From Lemma~\ref{Lem aN1} we therefore obtain that
\[\frac{\log \si_k(N)}{N} \le \frac{m \log a_{2^k}(1)}{N}
=\frac{m \log (v_k+1)+m\cdot 2^k \ga_k}{N}=\frac{\log(v_k+1)}{2^{k}}+\ga_k,\]
which proves the upper estimate for $\si_k(N)$. To get the lower estimate, let
$p_N(z)$ be an arbitrary $k$-fold symmetric Hurwitz stable polynomial of degree $N$ with
positive integer coefficients. We then have
\[
p_N(1)  =  p_{N/2}(v_1) = v_1^{N/4}p_{N/4}(v_2)=v_1^{N/4}v_2^{N/8} p_{N/8}(v_3)
\]
and so on, terminating with
\[\log p_N(1)=\frac{N}{4}\log v_1+\frac{N}{8}\log v_2 + \cdots + \frac{N}{2^k}\log v_k+
\log p_{N/2^k}(v_k),\]
which is the same as
$\log p_N(1)= N \ga_k + \log p_{N/2^k}(v_k)$.
From Theorem \ref{Theo 3.1} we now deduce that
\[\log p_N(1)\ge N \ga_k +\frac{N}{2^{k+1}}\log(v_k^2+1) \ge N\ga_k.\]
Taking the exponential we arrive at the asserted lower estimates. $\;\:\square$

\vsg
For $k=1,2,3$ the bounds provided by Theorem \ref{Theo kfold} read as follows.

\[\begin{array}{|c|c|}
\hline
k &  \Big((v_k^2+1)^{1/2^{k+1}}\Big) \ee^{\ga_k} \le \si_k(N)^{1/N} \le (v_k+1)^{1/2^k}\ee^{\ga_k}\\
\hline
  &    \\[-2ex]
1 & 1.4953\ldots = \sqrt[4]{5}  \le \si_k(N)^{1/N} \le \sqrt{3}=1.7320\ldots\\[1ex]
2 & 1.5233\ldots = \sqrt[8]{29} \le \si_k(N)^{1/N} \le \sqrt[4]{7}= 1.6265\ldots\\[1ex]
3 & 1.5340\ldots =\sqrt[16]{941} \le \si_k(N)^{1/N} \le \sqrt[8]{39}=1.5808\ldots\\
\hline
\end{array}\]

\vsg
Clearly, the $k$-fold sigmas $\si_k(N)$ also satisfy
the inequality $\si_k(N+M) \le \si_k(N)\si_k(M)$, and hence, by the argument of the proof of Theorem~\ref{Theo Fek},
the limits $\be_k$ of $\si_k(N)^{1/N}$ exist
as well.

\begin{cor} \label{Corbeta}
We have
\begin{eqnarray*}
& & 1.41 < 1.4142\ldots =\sqrt{2} \le \be \le \ee^\ga =1.5417\ldots < 1.55,\\
& & 1.49 < 1.4953\ldots =\sqrt[4]{5} \le \be_1 \le \ee^\ga =1.5417\ldots < 1.55,\\
& & 1.52 < 1.5233\ldots =\sqrt[8]{29} \le \be_2 \le \ee^\ga =1.5417\ldots < 1.55,\\
& & 1.53 < 1.5340\ldots = \sqrt[16]{941} \le \be_3 \le \ee^\ga =1.5417\ldots < 1.55.
\end{eqnarray*}
\end{cor}

{\em Proof.} The only thing we need to prove is the upper bound for $\si_k(N)$. So fix $k$
and take $N=2^{k+\ell}$ with $\ell=0,1,2, \ldots$. Then $\si_k(N) \le a_{2^{k+\ell}}(1)$,
and Lemma~\ref{Lem aN1} tells us that, given any $\eps >0$,
\[\frac{\log \si_k(N)}{N}\le \frac{\log(2\sqrt{k+\ell}+1)}{2^{k+\ell}}+\ga < \eps +\ga\]
whenever $\ell$ is large enough. This shows that the limit of $\si_k(N)^{1/N}$ does not
exceed $e^\ga$, as desired. $\;\:\square$

\begin{thm} \label{Theo 5.2}
If $N$ is divisible by $20$, then
\[\frac{1.41^N}{N} < p_{\max}(N) < 1.56^N\left(\frac{0.68}{\sqrt{N}}+0.97^N\right).\]
\end{thm}

{\em Proof.} To show the lower bound, let $p_N(z)$ be $c$-optimal. Then the sum of
the coefficients of $p_N(z)$ is at most $N p_{\max}$, and since this sum is greater
than $1.41^N$ by Theorem~\ref{Theo 5.1}, we conclude that $p_{\max} > 1.41^N/N$.

\vsk
The upper bound will follow once we have shown that the largest coefficient of
the polynomial $c_{20 k}(z):=c_{20}(z)^k$ is smaller than this bound, where $k=N/20$ and $c_{20}(z)$
is the polynomial~(\ref{c20}). The polynomial $c_{20}(z)$ is symmetric and hence
\[c_{20}(z)=z^{10}\left(c_{10}+c_9(z+z^{-1})+ \cdots + c_0(z^{10}+z^{-10})\right).\]
Thus,
\[c_{20 k}(z)=z^{10 k}\left(c_0(k)+c_1(k)(z+z^{-1})+ \cdots + c_{10 k}(k)(z^{10 k}+z^{-10 k})\right),\]
and the numbers $c_j(k)$ are the coefficients of $c_{20k}(z)$. For real $x$, we define
\[f(x)=c_{10}+c_9(\\e^{\ii x}+\ee^{-\ii x})+ \cdots + c_0 (\ee^{ 10 \ii x}+\ee^{-10 \ii x}).\]
Then $c_j(k)$ is just the $j$th Fourier coefficient of $f(x)^k$, that is,
\[c_j(k)=\frac{1}{2\pi}\int_{-\pi}^\pi f(x)^k \ee^{-\\i j x}dx.\]
It follows that
\[c_j(k)=|c_j(k)| \le \frac{1}{2\pi}\int_{-\pi}^\pi |f(x)|^k dx.\]
(Actually $f(x)>0$ for all $x$, but we don't need this.)
We have $f(0)=7167=:\si$ and the function
\[\frac{f(x)}{\si}=\frac{c_{10}}{\si}+\frac{2 c_9}{\si} \cos (x) + \cdots + \frac{2 c_0}{\si} \cos (10 x)\]
satisfies $0<f(x)/\si \le \ee^{-3.5 x^2}$ for $|x| <1$ and $|f(x)/\si| < 1/2$ for $1 \le |x| <\pi$.
Figure~5 shows the graphs of $f(x)/\si$ and $\ee^{-3.5 x^2}$.

\begin{figure}[thb] \label{Glock}
  \begin{center}
    \includegraphics[width=12cm]{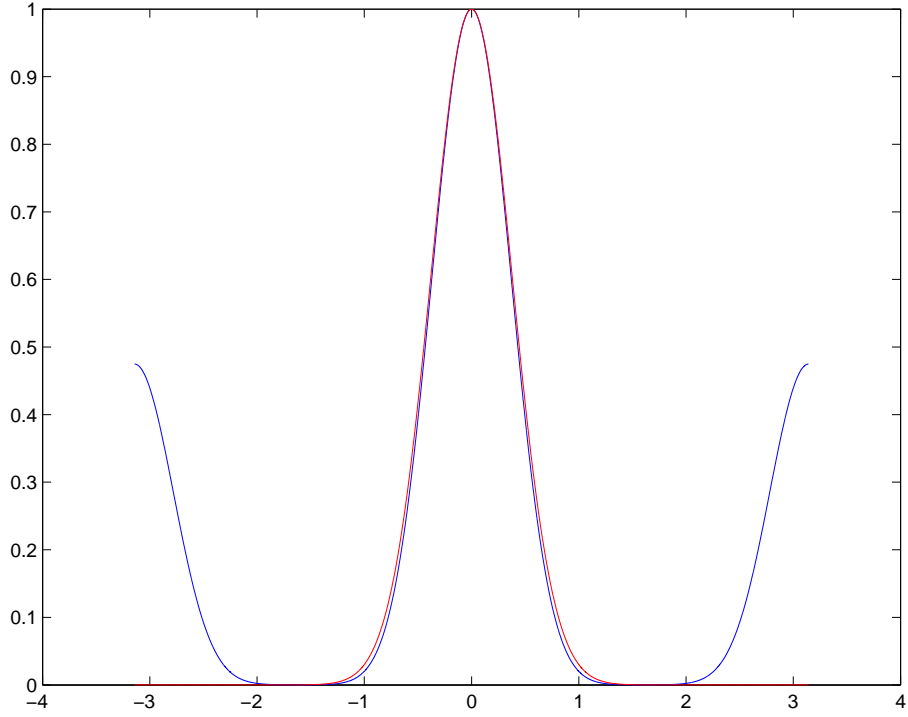}
    \caption{Graphs of $f(x)/\si$ (blue) and $\ee^{-3.5 x^2}$ (red).}
\end{center}
\end{figure}

\vsk
Consequently,
\begin{eqnarray*}
& & \frac{1}{2\pi \si^k}\int_{-\pi}^\pi f(x)^k dx < \frac{1}{2\pi}\int_{|x|<1}\ee^{-3.5 x^2} dx+\frac{1}{2\pi}\int_{1< |x| <\pi} 0.5^k dx\\
& & < \frac{1}{2\pi} \int_{-\iy}^\iy \ee^{-3.5 x^2} dx+ 0.5^k
= \frac{1}{\sqrt{14 \pi k}}+0.2^k.
\end{eqnarray*}
In summary we have
\[p_{\max}(N) \le \si^k\left(\frac{1}{\sqrt{14 \pi k}}+0.2^k\right).\]
Inserting $k=N/20$ and $\si=7167$ and taking into account that $7167^{1/20}< 1.56$ and $0.5^{1/20}<0.97$,
we arrive at the asserted upper bound. $\;\:\square$

\begin{rem} \label{Rem c20}
{\rm It is easy to find the asymptotics of the largest coefficient of the polynomials used in the preceding proof.}

\vsk
Let $c_{20}(z)$ be the polynomial~(\ref{c20}) and
put
\[\si=c_{20}(1)=7167, \quad \tau=\frac{26313}{7167}=3.6714 \ldots.\]
Then the maximum of the coefficients of $c_{20}(z)^k$ is
\[\frac{\si^k}{\sqrt{4\pi\tau k}}\,(1+o(1)).\]
\end{rem}

Indeed, we observed that the maximum
in question is
\[
\frac{\si^k}{2\pi}\int_{-\pi}^\pi \left(\frac{f(x)}{\si}\right)^k dx=\frac{\si^k}{2\pi}\int_{-\pi}^\pi \ee^{k g(x)} dx\]
with
\begin{eqnarray*}
\frac{f(x)}{\si} & = & \frac{c_{10}}{\si}+\frac{2 c_9}{\si} \cos (x) + \cdots + \frac{2 c_0}{\si} \cos (10 x)\\
& = & 1-\frac{1}{\si}(c_9+2^2 c_8+3^2 c_7 + \cdots + 10^2 c_0) x^2 +O(x^4) =1 - \tau x^2 +O(x^4).
\end{eqnarray*}
and
$g(x):=\log (f(x)/\si)= -\tau x^2+O(x^4)$. The function $g(x)$ is twice differentiable and attains it
maximum on $[-\pi,\pi]$ only at $x=0$. A well known theorem by Laplace therefore implies that
\[
\frac{\si^k}{2\pi}\int_{-\pi}^\pi \ee^{k g(x)} dx = \frac{\ee^{k g(0)}\si^k}{\sqrt{2\pi k |g''(0)|}} \,(1+o(1))
= \frac{\si^k}{\sqrt{2\pi k \cdot 2\tau}} \,(1+o(1)).
\]

Writing $k=N/20$ we get for the maximal coefficient of $c_N(z)=c_{20}(z)^k$ the asymptotics
\[(\sqrt[20]{7167})^N \,\sqrt{\frac{5}{\pi \tau N}}\,(1+o(1))\approx 1.5587 ^N\,\frac{0.6584}{\sqrt{N}}\,(1+o(1)), \]
which is in accordance with Theorem \ref{Theo 5.2}. The number $3.5$ we used in the proof of Theorem~\ref{Theo 5.2}
comes from the estimate $3.5 < \tau= 3.6714 \ldots$. $\;\:\square$

\bigskip
Albrecht B\"ottcher

Fakult\"at f\"ur Mathematik

TU Chemnitz

09107 Chemnitz

Germany

\bigskip
{\tt aboettch@mathematik.tu-chemnitz.de}


\begin{thebibliography}{99}

\bibitem{BiniEff}
D. A. Bini, G. Fiorentino, L. Gemignani, and B. Meini,
{\em Effective fast algorithms for polynomial spectral factorization}.
Numer. Algorithms 34 (2003), 217--227.

\bibitem{Fek}
M. Fekete,
{\em \"Uber die Verteilung der Wurzeln bei gewissen algebraischen Glei\-chungen mit
ganzzahligen Koeffizienten}.
Math. Zeitschrift 17 (1923), 228--249.

\bibitem{Law}
G. F. Lawler and L. N. Coyle,
{\em Lectures on Contemporary Probability}.
Amer. Math. Soc., Providence, RI, 1999.

\bibitem{Schur}
I. Schur,
{\em \"Uber die Verteilung der Wurzeln bei gewissen algebraischen Gleichungen mit
ganzzahligen Koeffizienten}.
Math. Zeitschrift 1 (1918), 377--402.

\bibitem{Trev}
V. Trevisan,
{\em Recognition of Hurwitz polynomials}.
SIGSAM Bull. 24 (1990), 26--32.
\end{thebibliography}
\end{document}